\documentclass{article}
\usepackage{amsmath,amssymb,small}
\usepackage[dvips]{graphicx}

\author{David R. E. Williams\footnote{The author was supported by EPSRC grant GRL67899.}}
\date{June 1999}
\title{Path-wise solutions of SDEs driven by L\'evy processes}

\begin{document}
\maketitle
\def\thefootnote{\fnsymbol{footnote}}
\def\@makefnmark{\hbox to\z@{$\m@th^{\@thefnmark}$\hss}}
\footnotesize\rm\noindent
\hspace*{16pt}\strut\footnote[0]{{\it AMS\/\ {\rm 1991} subject
classifications}.  60H20, 60G17, 60H05}\footnote[0]{{\it Key words and phrases}. L\'evy process, path integral, $p$-variation, area process, stochastic differential equations}
\normalsize\rm

\begin{abstract}
In this paper we show that a path-wise solution to the following integral equation 
$$
Y_t = \int_0^t f(Y_t)\;dX_t \qquad Y_0=a \in \R^d
$$
exists under the assumption that $X_t$ is a L\'evy process of finite $p$-variation for some $p \geq1$ and that $f$ is an $\alpha$-Lipschitz function for some $\alpha>p$.  There are two types of solution, determined by the solution's behaviour at jump times of the process $X$, one we call geometric the other forward.  The geometric solution is obtained by adding fictitious time and solving an associated integral equation.  The forward solution is derived from the geometric solution by correcting the solution's jump behaviour. 

L\'evy processes, generally, have unbounded variation.  So we must use a pathwise integral different from the Lebesgue-Stieltjes integral.  When $X$ has finite $p$-variation almost surely for $p<2$ we use Young's integral.  This is defined whenever $f$ and $g$ have finite $p$ and $q$-variation for $1/p+1/q>1$ (and they have no common discontinuities).  When $p>2$ we use the integral of Lyons.  In order to use this integral we construct the L\'evy area of the L\'evy process and show that it has finite $(p/2)$-variation almost surely. 
\end{abstract}

\section*{Introduction}\label{s.intro}

In this paper I give a path-wise method for solving the following integral equation: 
\begin{equation}\label{e.1}
Y_t=Y_0 + \int_0^t f(Y_t) \;dX_t\qquad Y_0=a\in\R^d.
\end{equation}
when the driving process is a L\'evy process.

Typically, a L\'evy process a.s. has unbounded variation.  The integral does not exist in a Lebesgue-Stieltjes sense.  However, the integral still makes sense as a random variable due to the stochastic calculus of semi-martingales developed by the Strasbourg school \cite{myrintsto}.  

The semi-martingale integration theory is not complete though.  There are processes of interest which do not fit into the semi-martingale framework, for example the fractional Brownian motion.  An alternative integral is provided by the path-wise approach studied by Lyons \cite{tel1}, \cite{tel5} and Dudley \cite{dud}.  The basis of their papers is that of Young \cite{you}, who showed that the integral 
\begin{equation}\label{e.youdef}
\int_0^t f\;dg
\end{equation}
is defined whenever $f$ and $g$ have finite $p$ and $q$-variation for $1/p+1/q>1$ (and they have no common discontinuities).  For a comprehensive overview of the theory we recommend the lecture notes of Dudley and Norvai{\v s}a \cite{dudnor}.  

Recently in \cite{miknor}, a system of linear Riemann-Stieltjes integral equations is solved when the integrator has finite $p$-variation for some $0<p<2$.  These results are contained in Theorem \ref{thm.1} where we allow non-linearity of the vector field $f$.  This is because our approach is an extension of the method of \cite{tel1}, \cite{tel5}.

The approach that I follow distinguishes two cases.  The first is when the process has finite $p$-variation a.s., for some $p<2$.  We use the Young integral \cite{you}.  In \cite{tel1} \eqref{e.1} is solved when $X_t$ is a continuous path of finite $p$-variation for some $p<2$.  

The second case is when the process has finite $p$-variation a.s., for some $p>2$.  The Young integral is only defined when $f$ and $g$ have finite $p$ and $q$-variation for $1/p+1/q>1$.  So an iteration scheme on the space of paths with finite $p$-variation does not work.  However, Lyons defined an integral against a continuous function of $p$-variation for some $p>2$ \cite{tel5}.  The integral is developed in the space of geometric multiplicative functionals (described in Appendix \ref{app.hom}).  The key idea is that we enhance the path by adding an area function to it.  If there is sufficient control of the pair, path {\em and} area, then the integral is defined.  The canonical example in \cite{tel5} is Brownian motion.  The area process enhancing the Brownian motion is the L\'evy area \cite[Ch.7, Sect.55]{levy}.  I show that there is an area process of a L\'evy process which has finite $(p/2)$-variation a.s..  

In order to solve \eqref{e.1} for a discontinuous function I add fictitious time during which linear segments remove the discontinuities, creating a continuous path.  By solving for the continuous path and then removing the fictitious time we recover a solution for the discontinuous path.  This is called a geometric solution.  A second type of solution is derined from the geometric solution which we call the forward solution. 

The first section treats the case where the discontinuous driving path has finite $p$-variation for some $p<2$.  The second section treats the case where the path has finite $p$-variation for some $p>2$ only.  The main proofs of the second section are deferred to the third section.  In the appendix I prove the homeomorphic flow property for the solutions when the driving path is continuous.  This is used in proving that forward solutions can be recovered from geometric solutions.

\section{Discontinuous processes - $p<2$}\label{s.pl2}

In this section we extend the results of \cite{tel1} to allow the driving path of \eqref{e.1} to have discontinuities.  The results are applied to sample paths of some L\'evy processes, those that have finite $p$-variation a.s. for some $p<2$.  Throughout this section $p\in[1,2)$ unless otherwise stated.

First, we determine the solution's behaviour when the integrator jumps.  There are two possibilities to consider:  the first is an extension of the Lebesgue-Stieltjes integral; the second is based on a geometric approach.

Suppose that the discontinuous integrator has bounded variation.  The solution $y$ would jump
$$
y_t-y_{t-} = f(y_{t-})\;(x_t-x_{t-})
$$
at a jump time $t$ of $x$.  If $x$ has finite $p$-variation for some $1<p<2$ we insert these jumps at the discontinuities of $x$.  We call a path $y$ with the above jump behaviour a forward solution.

The other jump behaviour we consider is the following:  When a jump of the integrator occurs we insert some fictitious time during which the jump is traversed by a linear segment, creating a continuous path on an extended time frame.  Then we solve the differential equation driven by the continuous path.  Finally we remove the fictitious time component of the solution path.  We call this a geometric solution because the solution has an 'instantaneous flow' along an integral curve at the jump times.  This jump behaviour has been considered before by \cite{marc} and \cite{kur2}.  

The disadvantage of the first approach is that the solution does not, generally, generate a flow of diffeomorphisms \cite{lea2}.

In this section we prove the following theorem:

\begin{thm}\label{thm.1}
Let $x_t$ be a discontinuous function of finite $p$-variation for some $p<2$.  Let $f$ be an $\alpha$-Lipschitz vector field for some $\alpha>p$.  Then there exists a unique geometric solution to the integral equation
\begin{equation}\label{e.2}
y_t = y_0 + \int_0^t f(y_t)\;dx_t \qquad y_0=a \in\R^d.
\end{equation}
With the above assumptions, there exists a unique forward solution as well.
\end{thm}

Before proving the theorem we recall the definitions of $p$-variation and $\alpha$-Lipschitz:

\begin{defn}
The $p$-variation of a function ${ x }(s)$ over the interval $[0,t]$ is defined as follows:
\begin{equation*}
{\Vert}{ x }{\Vert}_{_{p,[0,t]}}=\bigg\{\sup_{\pi\in\pi[0,t]}\quad{\sum_\pi}{\vert}{ x }(t_k)-{ x }(t_{k-1}){\vert}^p\bigg\}^{1\over
p}
\end{equation*}
where $\pi[0,t]$ is the collection of all finite partitions of the interval
$[0,t]$.  
\end{defn}
\begin{rem}
This is the strong $p$-variation.  Usually probabilists use the weaker form where the supremum is over partitions restricted by a mesh size which tends to zero.   
\end{rem}

\begin{defn}
  A function $f$ is in $\lip(\alpha)$ for some $\alpha>1$ if
\begin{equation*}
\norm{f}{\infty} < \infty \;\hbox{and}\; {{\partial f}\over{\partial x_j}} \in \hbox{Lip}(\alpha -1)\quad j=1,\dots ,d\,.
\end{equation*}
Its norm is given by
\begin{equation*}
\norm{f}{\lip(\alpha)}  \dfn  \norm{f}{\infty} + \sum_{j=1}^d \lnorm{{{\partial{f}}\over{\partial{x_j}}}}{\lip(\alpha -1)}	\qquad\qquad\hbox{for}\;\alpha>1.
\end{equation*}
\end{defn}

This is Stein's \cite{ste} definition of $\alpha$-Lipschitz continuity for $\alpha>1$.  It extends the classical definition:  $f$ is in $\lip(\alpha)$ for some $\alpha\in(0,1]$ if
\begin{equation*}
\snorm{f(x) - f(y)} \leq K \snorm{x-y}^{\alpha}
\end{equation*}
with norm
\begin{equation*}
\norm{f}{\infty} + \sup_{x \neq y} \frac{\snorm{f(x)-f(y)}}{\snorm{x-y}^{\alpha}}.
\end{equation*}

\subsection{Geometric Solutions.}
In this subsection we define a parametrisation for a c\`adl\`ag path $x$ of finite $p$-variation.  The parametrisation adds fictitious time allowing the traversal of the discontinuities of the path $x$.  We prove that the resulting continuous path $x^{\delta}$ has the same $p$-variation that $x$ has.  We solve \eqref{e.2} driven by $x^{\delta}$ using the method of Lyons \cite{tel1}.  Then we get a geometric solution of \eqref{e.2} by removing the fictitious time (i.e. by undoing the parametrisation).
 
\begin{figure}
\begin{center}
\includegraphics[angle=270,width=7cm]{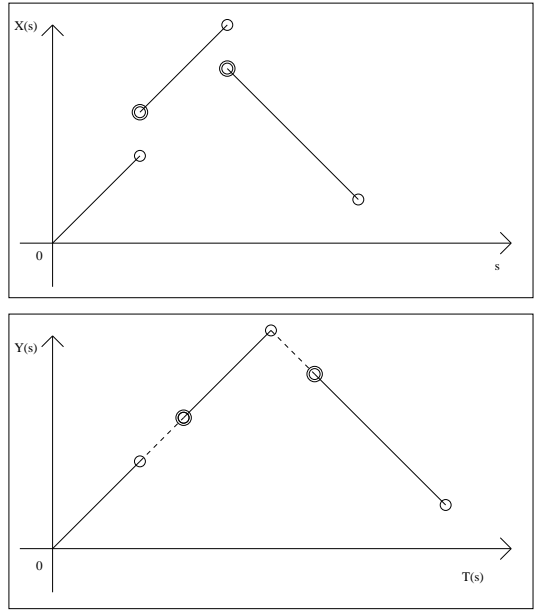}
\caption[Parametrisation of a discontinuous path.]{Let $x_s$ be a discontinuous path of bounded variation $(p=1)$.  Define a map $\tau(s)$ inserting fictitious time for the discontinuities of $x$.  Define a parametrisation $(t(s),y(s))$ in the manner of \eqref{e.xdel}.  $(t(s),y(s))$ traverses the jumps of $x$ during the fictitious time.}\label{fig.1}
\end{center}
\end{figure}

\begin{defn}\label{d.tau}
Let ${x}$ be a c\`adl\`ag path of finite $p$-variation.  Let $\delta>0\,,$ for each $n \geq 1\,,$ let $t_n$ be the time of the $n$'th largest jump of ${x}$.  We define a map $\tau^{\delta} : [0,T]\rightarrow [0, T+\delta
\sum_{i=1}^{\infty}\vert j(t_i) \vert^p] $ (where $j(u)$ denotes
the jump of the path ${x}$ at time u) in the following way:
\begin{equation}\label{e.tau}
\tau^{\delta}(t) = t + \delta \sum_{n=1}^{\infty} \vert j(t_n)\vert^p \chi_{\{t_n \leq t\}}(t).
\end{equation}
The map $\tau^{\delta} : [0,T]\rightarrow[0,\tau^{\delta}(T)]$ extends the time interval into one on where we define the continuous process ${x}^{\delta}(s)$
\begin{align}\label{e.xdel}
{x}^{\delta}(s)&\nonumber\\
=& \left\{
\begin{array}{ll}
 {x}(t) &\text{if}\; s = \tau^{\delta}(t) ,\\
 {x}(t_n^-) + (s - \tau^{\delta}(t_n^-)) j(t_n)\delta^{-1}\vert j(t_n) \vert ^{-p} &\text{if}\; s\in[\tau^{\delta}(t_n^-) \tau^{\delta}(t_n)) .
\end{array}
\right.
\end{align}
\end{defn}
\begin{rems}
$$ $$
\vspace{-10mm}
\begin{enumerate}
\item $(s,\,{x}^{\delta}_s),\;s\in [0,\tau^{\delta}(T)]$ is a parametrisation of the driving path ${x}$. 
\item The terms $\snorm{j(t_n)}^{p}$ in \eqref{e.tau} ensure that the addition of the fictitious time does not make $\tau^{\delta}(t)$ explode. 
\item In Figure \ref{fig.1} we see an example of a parametrisation of a discontinuous path $x_s$ in terms of the pair $(t(s),\,y(s))$. 
\end{enumerate} 
\end{rems}
The next proposition shows that the above parametrisation has the same $p$-variation as the original path, on the extended time frame $[0,\tau^{\delta}(T)]$.
\begin{propn}\label{p.pvarn}
Let ${x}$ be a c\`adl\`ag path of finite $p$-variation.  Let ${x}^{\delta}$ be a parametrisation of ${{x}}$ as above.  Then 
\begin{equation*} 
\variation{{x}^{\delta}}{p}{\tau^{\delta}(T)} =
\variation{{x}}{p}{T} \qquad \forall \delta > 0 .
\end{equation*}
\end{propn}
\proof
Let $\pi_0$ be a partition of $[0,\tau^{\delta}(T)]$.  Let
\begin{equation*}
V_{x^{\delta}}(\pi_0) = \sum_{\pi_0}\vert {x}^{\delta}(t_i) - {x}^{\delta}(t_{i-1}) \vert^p
\end{equation*}
We show that we increase the value of $V_p(\pi_0)$ by moving points lying on the jump segments to the endpoints of those segments.

Let $t_{i-1},t_i,t_{i+1}$ be three neighbouring points in the partition $\pi_0$ such that $t_i$ lies in a jump segment.  Consider the following term:
\begin{equation}\label{e.dom}
\snorm{x^{\delta}_{t_i} - x^{\delta}_{t_{i-1}}}^p + \snorm{x^{\delta}_{t_{i+1}}-x^{\delta}_{t_i} }^p.
\end{equation}
We show that \eqref{e.dom} is dominated by replacing $x^{\delta}_{t_i}$ by one
of $x^{\delta}_l$ and $x^{\delta}_r$, where $l$ and $r$ denote the left and right endpoint of the jump segment containing $t_i$.  

For simplicity we set $a=x^{\delta}_{t_{i-1}}, b=x^{\delta}_{t_{i+1}}$ and $c=x^{\delta}_l$.  Let
\begin{equation*}
L \dfn \{{c}+k{x} : k\in(0,1),\quad {c},{x}\in\R^d,{x}\neq0\} , {a},{b}\in\R^d\backslash L .
\end{equation*}
Let the function $f:[0,1]\rightarrow (0,\infty)$ be defined by
\begin{equation*}
f(k) = \vert {a} - {d}\vert^p + \vert {d} - {b}\vert^p \qquad, {d} = {c}+k{x}.
\end{equation*}
Then $f\in C^2[0,1]$ and one can show that $f'' \geq 0$ on $(0,1)$ when $p
\geq 1$.  To conclude the proof we move along the partition replacing $t_i$ which lie in the jump segments by new points $t_i^{\prime}$ that increase $V_{x^{\delta}}(\pi_0)$. The
partition $\pi_0$ is replaced by a partition $\pi_0'$ whose points lie on the pre-image of $[0,\tau^{\delta}(T)]$.  Therefore we have
\begin{equation*}
V_{x^{\delta}}(\pi_0) \leq V_{x^{\delta}}(\pi_0') = V_{x}(\pi_0').
\end{equation*}
Hence $\variation{{x}^{\delta}}{p}{\tau^{\delta}(T)} = \variation{{x}}{p}{T}$.\endproof

\begin{thm}\label{thm.geo}
Let ${x}$ be a c\`adl\`ag path with finite $p$-variation for some $p<2$.  Let $f$ be a $\lip(\gamma)$ vector field on $\R^n$ for some $\gamma> p$. Then there exists a unique geometric solution ${y}$, having finite $p$-variation which solves the differential equation
\begin{equation}\label{e.21}
d{y}_t = f({y}_t)\; d{x}_t \qquad {y}_0 = {a} \in \R^n.
\end{equation}
\end{thm}
\proof
Let $x^{\delta}$ be the parametrisation given in \eqref{e.xdel}.  The theorem of section three of \cite{tel1} proves that there is a continuous solution ${y}^{\delta}$ which solves \eqref{e.2} on $[0,\tau^{\delta}(T)]$.  Then $(s,\,{y}^{\delta}_s)$ is a parametrisation of a c\`adl\`ag path ${y}$ on  $[0,\,T]$.  

The solution is well-defined.  To see this, consider two parametrisations of ${x}$ and note that there exists a monotonically increasing function $\lambda_s$ such that 
\begin{equation*}
(s,{x}^{\delta}_s) = (\lambda_s,{x}^{\nu}_{\lambda_s}).\eendproof
\end{equation*}

\subsection{Forward Solutions.}\label{ss.for}
In this subsection we show how to recover forward solutions from geometric solutions.  The idea behind our approach is to correct the jump behaviour of the geometric solution using a Taylor series expansion Lemma \ref{l.tay}.  The correction terms are controlled by
$$
\sum_{i=1}^{\infty} \snorm{x_{t_i}-x_{t^-_i}}^2
$$     
which is finite due to the finite $p$-variation of the path $x$.

In the case where the driving path has only a finite number of jumps we note that the forward solution can be recovered trivially.  It is enough to mark the jump times of $x$ and solve the differential equation on the components where $x$ is continuous, inserting the forward jump behaviour when the jumps occur.  It remains to show that the forward solution exists when the driving path has a countably infinite number of jumps.  The method we use requires the following property of the geometric solution:

\begin{thm}\label{thm.hom}
Let $x$ be a continuous path of finite $p$-variation for some $p>1$.  Let $f$ be in $\lip(\alpha)$ for some $\alpha>p$.  The maps $(\pi_t)_{t\geq0}:\R^n\rightarrow \R^n$ obtained by varying the initial condition of the following differential equation generate a flow of homeomorphisms:
\begin{equation}\label{eq.fde}
d\pi_t = f(\pi_t)\;dx_t \qquad \pi_0 = Id, \quad (\text{the identity map}).
\end{equation}
\end{thm}

We leave the proof of Theorem \ref{thm.hom} until Appendix \ref{app.hom}.  We note the uniform estimate
\begin{equation}\label{e.simhom}
\sup_{0\leq t \leq T}\slnorm{\pi^{a}_t-\pi^{b}_t} \leq C(T) \slnorm{a-b}. 
\end{equation}

The following lemma will enable estimates to be made when the geometric jumps are replaced by the forward jumps:

\begin{lem}\label{l.tay}
Let $x$ be a c\`adl\`ag path with finite $p$-variation.  Let $f$ be in $\lip(\alpha)$ for some $\alpha>p$.  Let $\Delta y_i $ (resp. $\Delta z_i $) denote the geometric (resp. forward) solution's jump which correspond to $\Delta x_i$, the $i$'th largest jump of $x$.  Then we have the following estimate on the difference of the two jumps:
\begin{equation*}
\norm{\Delta y_i - \Delta z_i}{\infty} \leq K \vert \Delta x_i \vert^2
\end{equation*}
where the constant $K$ depends on $\norm{f}{\lip(\alpha)}$.
\end{lem}

\proof
Parametrise the path $x$ so that it traverses its discontinuity in unit time.  Solve geometrically over this interval with the solution having initial point $a$.  Note that the forward jump is the first order Taylor approximation to the geometric jump.  Then
\begin{eqnarray}
y_1(a) = y_0(a) &+& \left.\deriv{y_s(a)}{s}\right\vert_{s=0} + \half\;\left.\deriv{^2y_s(a)}{s^2}\right\vert_{s=\theta}\quad \hbox{for some}\;0<\theta<1\nonumber\\
&=& z_1(a) + \half\; \left.\deriv{^2y_s(a)}{s^2}\right\vert_{s=\theta}.
\end{eqnarray}
We estimate the second order term by
\begin{eqnarray}
\llnorm{\half\;\deriv{^2y_s(a)}{s^2}}{\infty} &=& \llnorm{\half\;\frac{d}{ds}f(y_s(a))(\Delta x_i)}{\infty}\nonumber\\
&\leq& \half\;\lnorm{\nabla f}{\infty} \lnorm{f}{\infty} \vert  \Delta x_i\vert^2\nonumber\\
&\leq& \half\; \lnorm{f}{\lip(\alpha)}^2 \vert  \Delta x_i\vert^2
\end{eqnarray}
Both $\norm{\nabla f}{\infty}$ and $\norm{f}{\infty}$ are finite because $f$ is $\lip(\alpha)$ for some $\alpha >p \geq 1$.\endproof

\begin{thm}\label{thm.24}
Let $x$ be a c\`adl\`ag path with finite $p$-variation.  Let $f$ be in $\lip(\alpha)$ for some $\alpha >p$.  Then there exists a unique forward solution to the following differential equation:
\begin{equation}\label{e.unifor}
dz_t = f(z_t)\;dx_t     \qquad\qquad z_0=a.
\end{equation}
\end{thm}

\proof
By Theorem \ref{thm.hom} there exists a unique homeomorphism $y$ which solves 
\begin{equation*}
dy_t = f(y_t)\;dx_t     \qquad\qquad y_0=a
\end{equation*}
in a geometric sense.

Label the jumps of $x$ by $j_{x}=\{j_i\}_{i=1}^{\infty}$ according to their decreasing size.  Let $z^n$ denote the path made by replacing the geometric jumps of $y$ corresponding to $\{j_i\}_{i=1}^{n}$ by the forward jumps $\{f()\,(\Delta x_i)\}_{i=1}^n$.  We show that the $(z^n)_{n\geq1}$ have a uniform limit.

We order the corrected jumps chronologically, say $\{t_i\}_{i=1}^n$.  Then we estimate the following term using Lemma \ref{l.tay} and the uniform bound on the growth of $y$ given in \eqref{e.simhom}:
\begin{eqnarray}
\snorm{z^n_s(a) - y_s(a)} &\leq& \sum_{i=1}^n \snorm{y_{t_i,s}(z^n_{t_i}(a)) - y_{t_i,s}(y_{t_{i-1},t_i}(z^n_{t_{i-1}}(a)))}\nonumber\\
&\leq& C(T) \sum_{i=1}^n \snorm{z^n_{t_i}(a) - y_{t_{i-1},t_i}(z^n_{t_{i-1}}(a))}\nonumber\\
&\leq& C^2(T)\; K\;\sum_{i=1}^{\infty} \snorm{\Delta x_i}^2.
\end{eqnarray}
So we have the uniform estimate
\begin{equation}\label{e.234}
\norm{z^n-y}{\infty} \leq K(C_3(T),\norm{f}{\lip(\alpha)}) \sum_{i=1}^{\infty} \vert \Delta x_i \vert^2 < \infty        \qquad \forall n \geq 1.
\end{equation}

We use an analogous bound to get Cauchy convergence of $\{z^n\}_{n\geq1}$.  Let $m, r \geq 1$.
\begin{equation*}
\norm{z^m - z^{m+r}}{\infty} \leq K (C(T,z^m),\norm{f}{\lip(\alpha)}) \sum_{i=m+1}^{\infty} \vert \Delta x_i \vert^2\,.
\end{equation*}
One notes that $\{C(T,z^m)\}$ are uniformly bounded, because of the boundedness of $C(T)= C(T,y)$ and the Lipschitz condition on $f$.  Therefore we have the following estimate:
\begin{equation*}
\norm{z^m - z^{m+r}}{\infty} \leq L \sum_{i=m+1}^{\infty} \snorm{\Delta x_i}^2.
\end{equation*}
This implies that $\{z^n\}$ are Cauchy in the supremum norm because $x$ has finite $p$-variation $(p<2)$ which implies that $\sum_{m+1}^{\infty} \snorm{\Delta x_i}^2$ tends to zero as $m$ increases.\endproof

\begin{rem}
Theorems \ref{thm.24} and \ref{thm.geo} combine to prove Theorem \ref{thm.1}.
\end{rem}

\begin{cor}\label{cor.p}
With the above notation, $z$ has finite $p$-variation.
\end{cor}
\proof
Let $s<t \in [0,T]$.
\begin{align*}
\snorm{z_t-z_s} &\leq \snorm{(z_t-z_s)-(y_t-y_s)} + \snorm{y_t-y_s}\\
\intertext{where $(y_t-y_s)$ is the increment of the geometric solution starting from $z_s$ driven by the path $x_t$ on the interval $[s,T]$.  Then}
\snorm{(z_t-z_s)-(y_t-y_s)} \leq& \quad C\;\sum_{\substack{ j_x\arrowvert_{[s,t]} }} \snorm{\Delta x_i}^2 \qquad\text{ and } \qquad\snorm{y_t-y_s} \leq \variationt{x}{p}{s}{t}\,,\\
\intertext{which implies that}
\snorm{z_t-z_s}^p &\leq 2^{p-1} \bigg\{ C^p \big(\sum_{\substack{j_x\arrowvert_{[s,t]}  }} \snorm{\Delta x_i}^2\big)^p + \variationt{x}{p}{s}{t}^p\bigg\}\,,\\
\intertext{hence}
\variation{z}{p}{T} &\leq 2^{(p-1)/p} \bigg\{ C^p \big(\sum_{\substack{ j_x\arrowvert_{[0,T]} }} \snorm{\Delta x_i}^2\big)^p + \variation{x}{p}{T}^p \bigg\}^{1/p} < \infty.
\end{align*}
\hspace{8cm}\endproof

\subsection{$p$-variation of L\'evy processes}\label{ss.plevy}

In this subsection we apply Theorem \ref{thm.1} to L\'evy processes which have finite $p$-variation a.s.. 

L\'evy processes are the class of processes with stationary, independent increments which are continuous in probability.  The class includes Brownian motion, although this process is atypical due to its continuous sample paths.  Typically a L\'evy process will be a combination of a deterministic drift, a Gaussian process and a jump process.  For further information on L\'evy processes we direct the reader to \cite{bertbk}.

The regularity of the sample paths of a L\'evy process has been studied intensively.  In the 1960's several people worked on the question of characterising the sample path $p$-variation.  The following theorem, due to Monroe, gives the characterisation:

\begin{thm}\cite[Theorem 2]{mon}\label{thm.monroe}
Let $(X_t)_{t\geq0}$ be a L\'evy process in $\R^n$ without a Gaussian part.  Let $\nu$ be the L\'evy measure.  Let $\beta$ denote the index of $X_t$, that is
\begin{equation}\label{e.index}
\beta \dfn \inf\bigg\{ \alpha >0 \,:\, \int_{\vert {y}\vert \leq 1} \vert {y}\vert^{\alpha}\, \nu(d{y}) < \infty\bigg\}
\end{equation}
and suppose that $1\leq \beta\leq2$.  If $ \gamma> \beta$ then 
\begin{equation}
\prob{\norm{X}{\gamma} < \infty} = 1 
\end{equation}
where the $\gamma$-variation is considered over any compact interval.
\end{thm}

\begin{rem}
Note that all L\'evy processes with a Gaussian part only have finite $p$-variation for $p>2$.
\end{rem}

\begin{cor}\label{c.levy}
Let $(X_t)_{t\geq0}$ be a L\'evy process with index $\beta<2$ and no Gaussian part.  Let $f$ be a vector field in $\lip(\alpha)$ for some $\alpha>p$. Then, a.s., the following stochastic differential equation has a unique forward and a unique geometric solution:
$$
dY_t = f(Y_t)\;dX_t\qquad Y_0=a.
$$
\end{cor}

\proof
The corollary follows immediately from Theorems \ref{thm.monroe} and \ref{thm.1}.\endproof

\section{Discontinuous processes - $p>2$}\label{s.pg2}

The goal of this section is to extend (Corollary \ref{c.levy}) to let any L\'evy process be the integrator of \eqref{e.1}.   

One problem we have is that the Young integral is no longer useful because we use a Picard iteration scheme which fails condition \eqref{e.youdef} when $p>2$.  However, we can use the method from \cite{tel5}.  To define the integral we need to provide more information about the sample path.  We do this by defining an area process of the L\'evy process.  Then we prove that the enhanced process (path {\em and} area) has finite $p$-variation Definition \ref{d.enhp}.  

We parametrise the enhanced process in an analogous manner to \eqref{e.xdel} (adding fictitious time).  Then we solve \eqref{e.1} in a geometric sense using the method for continuous paths $(p>2)$ given in \cite{tel5}.  Finally, forward solutions are obtained by jump correction as before.

Before enhancing $(X_t)_{t\geq0}$ we give an example which shows that there exist L\'evy measures with index two.  So a L\'evy process does not need a Gaussian part to have, a.s., finite $p$-variation only for $p>2$.

\begin{eg}
{\allowdisplaybreaks
One can define the following measures on $\R$:
\begin{align*}
\nu_k\;(dx) &\dfn \snorm{x}^{-3+1/k}\;dx\qquad \snorm{x}\in((k+1)^{-3(k+1)},k^{-3k}] \dfn J_k\\
\eta_m\;(dx) &\dfn \sum_{k=1}^m \nu_k\;(dx\cap J_k\cap (-J_k)) .
\end{align*}
We show that $\eta \dfn \lim_{m\rightarrow\infty} \eta_m$ is a L\'evy measure.  The integrability condition 
\begin{equation}\label{e.intcon}
\int_{\snorm{x}\leq 1} \snorm{x}^2 \eta\;(dx) < \infty.
\end{equation}
must be satisfied.
\begin{align*}
\int_{\snorm{x}\leq 1} \snorm{x}^2 \eta_m\;(dx) &= 2 \int_0^1 \sum_{k=1}^m x^{-1+1/k} \chi_{J_k}(x)\;dx = 2 \sum_{k=1}^m \bigg[ k x^{1/k} \bigg]_{(k+1)^{-3(k+1)}}^{k^{-3k}}\\
&= 2 \sum_{k=1}^m k\bigg\{k^{-3} - (k+1)^{-3(1+1/k)}\bigg\}\\
&\leq 2 \sum_{k=1}^m k\bigg\{k^{-3} - 2^{-3(1+1/k)} k^{-3(1+1/k)}\bigg\}\\
&=  2 \sum_{k=1}^m k^{-2}\bigg\{1 - 2^{-3(1+1/k)} k^{-3/k}\bigg\}\\
&< C \sum_{k=1}^{\infty} k^{-2}\;<\infty, 
\end{align*}
where $C$ is some suitable constant.  We take the limit as $m$ tends to infinity on the left hand side to prove \eqref{e.intcon}.

Now we show that 
\begin{equation}\label{e.intpvar}
\int_{\snorm{x}\leq 1} \snorm{x}^{\alpha} \eta\;(dx) = \infty
\end{equation}
for all $\alpha<2$.}
{\allowdisplaybreaks
Fix $\alpha<2$.  Define the following number:
\begin{align*}
m(\alpha) &\dfn \inf\{ k\;:\; \alpha+1/k <2\} < \infty \qquad\text{as}\;\alpha<2.\\
\intertext{Let $m>m(\alpha)$.  Then}
&\int_{\snorm{x}\leq 1} \snorm{x}^{\alpha} \eta_m\;(dx) \\
&\geq 2 \sum_{k=m(\alpha)}^m {{1}\over{(\alpha+1/k-2)}}\bigg\{ k^{-3k(\alpha+1/k-2)}-(k+1)^{-3(k+1)(\alpha+1/k-2)}\bigg\}\\
&= 2 \sum_{k=m(\alpha)}^m {{1}\over{(2-(\alpha+1/k))}}\bigg\{ (k+1)^{-3(k+1)(\alpha+1/k-2)}- k^{-3k(\alpha+1/k-2)}\bigg\}\\
&\geq {{2}\over{2-\alpha}} \sum_{m(\alpha)}^m \bigg\{ (k+1)^{3(k+1)(2-(\alpha+1/k))}- k^{3k(2-(\alpha+1/k))}\bigg\}\\
&\rightarrow \infty \qquad\text{as}\; m\rightarrow \infty\,.
\end{align*}
}
This proves that the index $\beta$ of $\eta$ equals two.  Theorem \ref{thm.monroe} implies that the pure jump process associated to the L\'evy measure $\eta$ a.s. has finite $p$-variation for $p>2$ only.
\end{eg}  

The following theorem gives a construction of the L\'evy area of the L\'evy process $(X_t)_{t\geq0}$.  The L\'evy area process and the L\'evy process form the enhanced process which we need in order to use the method of Lyons \cite{tel5}.

\begin{thm}\label{thm.areaconv}
The $d$-dimensional L\'evy process $(X_t)_{t\geq0}$ has an anti-symmetric area process 
\begin{equation*}
(A_{s,t})^{ij} \dfn \half \int_s^t X^i_{u-} \circ dX^j_u - X^j_{u-} \circ dX^i_u \qquad i,j = 1,2. \qquad\as
\end{equation*}
\end{thm}
The proof is deferred to Section \ref{s.proofs}.

\begin{thm}\label{thm.parea}
The L\'evy area of the L\'evy process $(X_t)_{t\geq0}$ a.s. has finite $(p/2)$-variation for $p>2$.  That is
\begin{equation*}
\sup_{\pi} \big(\sum_{\pi} \snorm{A_{t_{k-1},t_k}}^{p/2}\big)^{2/p}\;<\infty\qquad\as
\end{equation*}
where the supremum is taken over all finite partitions $\pi$ of $[0,T]$.
\end{thm}

\noindent The proof is deferred to Section \ref{s.proofs}.

Now we parametrise the sample paths of $(X_t)_{t\geq0}$ as before \eqref{e.xdel}.  

\begin{propn}
Parametrising the process $(X_t)_{t\geq0}$ does not affect the area process' $(p/2)$-variation.
\end{propn}

\proof
The proof is similar to the proof of Proposition \ref{p.pvarn}.  One can show that if $\lambda$ lies in a jump segment then 
$$
\snorm{A_{s,\lambda}}^{(p/2)} +\snorm{A_{\lambda,t}}^{(p/2)}\qquad s<\lambda<t
$$
is maximised when $\lambda$ is moved to one of the endpoints of the jump segment.\endproof

With the parametrisation of the path and the area we can define the integral in the sense of Lyons \cite{tel5}.  Consequently we have the following theorem:

\begin{thm}
Let $(X_t)_{t\geq0}$ be a L\'evy process with finite $p$-variation for some $p>2$.  Let $f$ be in $\lip(\alpha)$ for some $\alpha>p$.  Then there exists, with probability one, a unique geometric and a unique forward solution to the following integral equation:
\begin{equation}\label{e.ass}
Y_t = Y_0 +\int_0^t f(Y_t)\;dX_t\qquad Y_0=a\in\R^d.
\end{equation}
\end{thm}

\begin{rem}
When constructing the forward solution it is necessary that the sum
$$
\sum_{n=1}^{\infty} \snorm{\Delta X_n}^2
$$
remains finite.  This is guaranteed by the requirement on L\'evy measures to satisfy
$$
\int_{\snorm{x}\leq1}\;\snorm{x}^2\wedge1\;\nu(dx) < \infty.
$$
\end{rem}

\section{Proofs of Theorem \ref{thm.areaconv} and Theorem \ref{thm.parea}}\label{s.proofs}

For clarity throughout this section we assume that the L\'evy process $(X_t)_{t\geq0}$ is two dimensional and takes the following form:
\begin{equation}\label{e.form}
X_t=B_t+\int_{\snorm{x}\leq1}\;x\;(N_t(dx)-t\nu(dx)).
\end{equation}
That is, $(X_t)_{t\geq0}$ is a Gaussian process with a compensated pure jump process, whose L\'evy measure is supported on $(x\in\R^2\, :\, \snorm{x}\leq1)$.

\begin{propn}\label{prop.areaconv}
The $d$-dimensional L\'evy process $(X_t)_{t\geq0}$ has an anti-symmetric area process
\begin{equation*}
(A_{s,t})^{ij} \dfn \half \int_s^t X^i_{u-} \circ dX^j_u - X^j_{u-} \circ dX^i_u \qquad i,j = 1,2. \qquad\as
\end{equation*}
For fixed $s<t$ we obtain the area process by the following limiting procedure:
\begin{equation*}
(A_ {s,t})^{ij} = \lim_{n\rightarrow \infty} \sum_{m=0}^n \sum _{\substack{ k=1,\\ odd }} ^{2^m -1} A^{i,j}_{k,m}\qquad \as
\end{equation*}
where $A^{ij}_{k,m}$ is the area of the $(ij)$-projected triangle with vertices
\begin{equation*}
X(u_{(k+1)/2,m-1}),\; X(u_{(k-1)/2,m-1}), \;X(u_{k,m})
\end{equation*}
where $u_{k,m} \dfn s+ k2^{-m}(t-s)$.  Also we have the second order moment estimate
\begin{equation}\label{e.2ndmom}
\expn{(A_{s,t}^{ij})^2} \leq C(\nu) (t-s)^2.
\end{equation}
\end{propn}

\proof
We define $A_{s,t}(n)$ 
\begin{align*}
A_{s,t}(n) &\dfn \half \sum_{k=0}^{2^n-1} (X^{(1)}(u_{k,n})-X^{(1)}(s))(X^{(2)}(u_{k+1,n})-X^{(2)}(u_{k,n}))\\
&\qquad \qquad\qquad- (X^{(2)}(u_{k,n})-X^{(2)}(s))(X^{(1)}(u_{k+1,n})-X^{(1)}(u_{k,n}))\\
& = \sum_{k=0}^{2^n-1} B_{k,n}
\end{align*}
where $B_{k,n}$ is the (signed) area of the triangle with vertices
\begin{equation*}
X(s),\; X(u_{k,n}),\; X(u_{k+1,n}).
\end{equation*}
By considering the difference between $A_{s,t}(n)$ and $A_{s,t}(n+1)$ we see that
\begin{equation*}
B_{2k,n+1} + B_{2k+1,n+1} - B_{k,n} 
\end{equation*}
is the area of the triangle with vertices 
\begin{equation*}
X(u_{k,n}),\; X(u_{k+1,n}),\; X(u_{2k+1,n+1})
\end{equation*}
which we denote by $A_{k,n}$.  We re-order $A_{s,t}(n)$
\begin{align*}
A_{s,t}(n) &= \half \sum_{m=0}^n \sum _{\substack{ k=1, \\ odd }} ^{2^m-1} \bigg(X(u_{k,m}) - d_{k,m}\bigg)\\
&\qquad\qquad\qquad\qquad\otimes\bigg(X(u_{(k+1)/2,m-1})-X(u_{(k-1)/2,m-1})\bigg)\\
&= \half  \sum_{m=0}^n \sum _{\substack{ k=1, \\ odd }} ^{2^m-1} A_{k,m}
\end{align*}
where $d_{k,m} \dfn 1/2\, (X(u_{(k+1)/2,m-1}) + X(u_{(k-1)/2,m-1}))$.  The convergence to the area process is completed using martingale methods.

Let $\sigf_n \dfn \sigma \big( X(u_{k,n}) \;:\; k=0,\ldots,2^n\big).$  Then
\begin{lem}
\begin{equation}\label{e.cent}
\condexp{X(u_{k,m})}{\sigf_{m-1}} = d_{k,m} \qquad \as
\end{equation}
\end{lem}
\proof
For ease of presentation we let 
\begin{align*}
U_1 &\dfn X(u_{k,m}) - X(u_{(k-1)/2,m-1})\\
U_2 &\dfn X((u_{(k+1)/2,m-1}) - X(u_{k,m}).
\end{align*}
Then 
\begin{align}
&\condexp{X(u_{k,m}) - d_{k,m}}{\sigf_{m-1}} \nonumber\\
&= \condexp{X(u_{k,m}) - d_{k,m}}{X(u_{(k-1)/2,m-1}),X(u_{(k+1)/2,m-1})}\nonumber\\
&= \half\; \condexp{U_1 - U_2}{X(u_{(k-1)/2,m-1}),X(u_{(k+1)/2,m-1})} \nonumber
\end{align}
{\allowdisplaybreaks

Using the stationarity and the independence of the increments of $X$ we see that $U_1$ and $U_2$ are exchangeable, that is
\begin{equation*}
\prob{U_1 \in A ,\; U_2 \in B} = \prob{U_2 \in A ,\; U_1 \in B}\qquad \forall A,B \in \borel(\R^2).
\end{equation*}
The exchangeability extends to the random variables
\begin{equation*}
\big(U_i \vert X(u_{(k-1)/2,m-1}),X(u_{(k+1)/2,m-1})\big)\qquad i=1,2.
\end{equation*}
We deduce that 
\begin{equation*}
\expn{U_1-U_2\;\vert\;X(u_{(k-1)/2,m-1}),X(u_{(k+1)/2,m-1})} = 0.\eendproof
\end{equation*}

Returning to the proof of Proposition \ref{prop.areaconv}, we compute the variance of $A_{k,m}$.  This will be used to show that 
\begin{equation*}
\sup_{n\geq 1}\, \expn{A_{s,t}(n)^2} < \infty.
\end{equation*}
\begin{align*}
\E\big(&A^2_{k,m}\big)\\
= \E\bigg(&\bigg([ X^{(1)}(u_{k,m}) - d^{(1)}_{k,m}][U_1^{(2)}+U_2^{(2)}] - [ X^{(2)}(u_{k,m}) - d^{(2)}_{k,m} ][U_1^{(1)}+U_2^{(1)}]\bigg)^2 \bigg)\\
&= \quart \expn{\bigg\{(U_1^{(1)}-U_2^{(1)})(U_1^{(2)}+U_2^{(2)})-(U_1^{(2)}-U_2^{(2)})(U_1^{(1)}+U_2^{(1)})\}^2}\\
&= \quart\expn{(U_1^{(1)}U_2^{(2)})^2 -2 U_1^{(1)}U_2^{(2)}U_2^{(1)}U_1^{(2)} +(U_2^{(1)}U_1^{(2)})^2 }\\
&\dfn \quad(1)\quad+\quad(2)\quad+\quad(3).
\end{align*}
We use the independence of the increments and It\^o's formula for discontinuous semi-martingales to compute $(1),(2)$ and $(3)$.
\begin{align*}
(1) &= \expn{(U_1^{(1)}U_2^{(2)})^2} = \expn{(U_1^{(1)})^2}\expn{(U_2^{(2)})^2}.
\end{align*}
By applying It\^o's formula and using the stationarity of the L\'evy process we find that
\begin{align*}
&(3)=(1)=2^{-2m} \;(t-s)^2 \int_{\snorm{x}\leq 1} \snorm{x_1}^2 \;\nu(dx)\;\int_{\snorm{x}\leq 1} \snorm{x_2}^2 \;\nu(dx).
\end{align*}
Another application of It\^o's formula gives
\begin{align*}
(2)&= -2 \expn{U_1^{(1)}U_2^{(2)}U_2^{(1)}U_1^{(2)}} = -2 \expn{U_1^{(1)}U_1^{(2)}}\expn{U_2^{(2)}U_2^{(1)}}\\
&= - 2^{-2m+1}\; (t-s)^2 \;\bigg(\int_{\snorm{x}\leq 1} x_1\,x_2 \;\nu(dx)\bigg)^2.
\end{align*}
Collecting the terms together we have the following expression:
\begin{align*}
&\expn{A^2_{k,m}}= C_0(\nu)\; 2^{-2m+1} \;(t-s)^2 \\
\intertext{where}
&C_0(\nu) \dfn \bigg\{\int_{\snorm{x}\leq 1} \snorm{x_1}^2 \;\nu(dx)\;\int_{\snorm{x}\leq 1} \snorm{x_2}^2 \;\nu(dx)-\bigg(\int_{\snorm{x}\leq 1} x_1\,x_2 \;\nu(dx)\bigg)^2\bigg\}.
\end{align*}

Now we estimate the following term:
\begin{align*}
\expn{A_{s,t}^2(n)} &= \expn{\bigg(\sum_{m=1}^n \sum _{\substack{ k=1,\\ odd }} ^{2^m -1} A_{k,m}\bigg)^2}\\
\intertext{which through conditioning and independence arguments equals}
&= \expn{\sum_{m=1}^n \sum _{\substack{ k=1,\\ odd }} ^{2^m -1} A_{k,m}^2}\\
&=  C_0(\nu) \sum_{m=1}^n \sum _{\substack{ k=1,\\ odd }} ^{2^m -1} 2^{-2m+1}(t-s)^2\\
&\leq C_0(\nu) \sum_{m=1}^{\infty} \sum _{\substack{ k=1,\\ odd }} ^{2^m -1} 2^{-2m+1}(t-s)^2 \dfn C(\nu) (t-s)^2.
\end{align*}
We use the martingale convergence theorem to deduce that a.s. there is a unique limit of $A_{s,t}(n)$.  Furthermore the last calculation implies that there is a moment estimate of the area process given by
\begin{equation*}
\expn{A^2_{s,t}} \leq  C(\nu) \;(t-s)^2.\eendproof
\end{equation*}

We note that there is another way that one could define an area process of a L\'evy process.  One could define the area process for the truncated L\'evy processes and look for a limit as the small (compensated) jumps are put in.  Using the above construction one can define $A^{\epsilon}_{s,t}$ for a fixed pair of times, corresponding to the L\'evy process $X^{\epsilon}$.  With the $\sigma$-fields $(\sigg^{\epsilon})_{\epsilon>0}$ defined by
$$
\gothic{G}^{\epsilon} \dfn \sigma ( X^{\delta} \,:\, \delta >\epsilon)\;\text{for} \;\epsilon >0
$$
we have the following proposition:

\begin{propn}\label{p.spamart}
$(A^{\epsilon}_{s,t})_{\epsilon>0}$ form a $(\sigg^{\epsilon})$-martingale.
\end{propn}
\proof
Let $\eta>\epsilon>0$.  By considering the construction of the area given above for the truncated processes $X^{\eta}$ and $X^{\epsilon}$ we look at the difference at the level of the triangles $A^{\eta}_{k,n}$ and $A^{\epsilon}_{k,n}$.  
\begin{align*}
\condexp{A^{\epsilon}_{k,n}-A^{\eta}_{k,n}}{\sigg^{\eta}} &= \E\bigg(A^{\eta,\epsilon}_{k,n}\\
&\quad+(X^{\eta,\epsilon}_{k,n} - d^{\eta,\epsilon}_{k,n})\otimes(X^{\eta}_{(k+1)/2,n-1}-X^{\eta}_{(k-1)/2,n-1})\\
&\quad+(X^{\eta}_{k,n} - d^{\eta}_{k,n})\otimes(X^{\eta,\epsilon}_{(k+1)/2,n-1}-X^{\eta,\epsilon}_{(k-1)/2,n-1})\big\vert \sigg^{\eta}\bigg)
\end{align*}
where the superscript $\eta,\epsilon$ signifies that the process is generated by the part of the L\'evy measure whose support is $(\epsilon,\eta]$.  Using the spatial independence of the underlying L\'evy process we have
\begin{align*}
&= \expn{A^{\eta,\epsilon}_{k,n}} + \expn{(X^{\eta,\epsilon}_{k,n} - d^{\eta,\epsilon}_{k,n})}\otimes(X^{\eta}_{(k+1)/2,n-1}-X^{\eta}_{(k-1)/2,n-1})\\
&\qquad\qquad+(X^{\eta}_{k,n} - d^{\eta}_{k,n})\otimes\expn{(X^{\eta,\epsilon}_{(k+1)/2,n-1}-X^{\eta,\epsilon}_{(k-1)/2,n-1})}\\
&= 0\,.\eendproof
\end{align*}

With the uniform control on the second moment of the martingale 
\begin{equation*}
\expn{(A^{\epsilon}_{s,t})^2} \leq  C(\nu) \;(t-s)^2 \qquad \forall \epsilon>0
\end{equation*}
we conclude that $A^{\epsilon}_{s,t}$ converges a.s. as $\epsilon\rightarrow0$.

The algebraic identity
\begin{equation}\label{e.algarea}
 A_{s,u} = A_{s,t} + A_{t,u} + \half [X_{s,t},X_{t,u}]  \qquad s<t<u
\end{equation}  
for the anti-symmetric area process $A$ generated by a piecewise smooth path $X$ extends to the area process of the L\'evy process.  This is due to \eqref{e.algarea} holding for the area processes $A^{\epsilon}$ of the truncated L\'evy processes $X^{\epsilon}$.

\begin{propn}\label{prop.parea}
The L\'evy area of the L\'evy process $(X_t)_{t\geq0}$ has finite $(p/2)$-variation for $p>2$ a.s..  That is
\begin{equation*}
\sup_{\pi} \big(\sum_{\pi} \snorm{A_{t_{k-1},t_k}}^{p/2}\big)^{2/p}\;<\infty\qquad\as
\end{equation*}
where the supremum is taken over all finite partitions $\pi$ of $[0,T]$.
\end{propn}
\proof
In Proposition \ref{prop.areaconv} we constructed the area process for a pair of times, a.s..  This can be extended to a countable collection of pairs of times, a.s..  In the proof below we assume that the area process has been defined for the times
\begin{equation*}
k2^{-n}T, (k+1)2^{-n}T \qquad k=0,1,\ldots,2^n-1, \qquad n\geq1.
\end{equation*}

The proof follows the method of estimation used in \cite{benterry}.  To estimate the area process for two arbitrary times $u<v$ we split up the interval $[u,v]$ in the following manner:

We select the largest dyadic interval $[(k-1)2^{-n}T,k2^{-n}]$ which is contained within $[u,v]$.  Then we add dyadic intervals to either side of the initial interval, which are chosen maximally with respect to inclusion in the interval $[u,v]$.  Continuing in this fashion we label the partition according to the lengths of the dyadics.  We note that there are at most two dyadics of the same length in the partition which we label $[l_{1,k},r_{1,k}]$ and $[l_{2,k},r_{2,k}]$ where $r_{1,k}\leq l_{2,k}$.  Then
\begin{equation*}
[u,v] = \bigcup_{k=1}^{\infty} \bigcup_{i=1,2} [l_{i,k},r_{i,k}].
\end{equation*}
We estimate $A_{u,v}$ using the algebraic formula \eqref{e.algarea}.
\begin{align*}
A_{l_{1,m},r_{2,m}} &= \sum_{k=1}^m \sum_{i=1,2} A_{l_{i,k},r_{i,k}}\\
& + \half \; \sum_{1\leq a\leq b \leq2} \sum_{1 \leq j < k \leq m} \bigg[ X_{r_{a,k}}-X_{l_{a,k}} , X_{r_{b,j}}-X_{l_{b,j}} \bigg].
\end{align*}
Noting that 
\begin{align*}
 &\sum_{1\leq a\leq b \leq2} \sum_{1 \leq j < k \leq m} \slnorm{\bigg[ X_{r_{a,k}}-X_{l_{a,k}} , X_{r_{b,j}}-X_{l_{b,j}}\bigg]} \\
&=\sum_{1\leq a\leq b \leq2} \sum_{1 \leq j < k \leq m} \slnorm{(X_{r_{a,k}}-X_{l_{a,k}})\otimes (X_{r_{b,j}}-X_{l_{b,j}})\\
&\qquad \qquad \qquad \qquad \qquad-\;(X_{r_{b,j}}-X_{l_{b,j}})\otimes (X_{r_{a,k}}-X_{l_{a,k}})}\\
&\leq \sum_{1\leq a\leq b \leq2} \sum_{1 \leq j < k \leq m}\slnorm{X_{r_{a,k}}-X_{l_{a,k}}}\slnorm{X_{r_{b,j}}-X_{l_{b,j}}}\\
&\leq \bigg( \sum_{k=1}^m \sum_{i=1,2} \slnorm{X_{r_{i,k}}-X_{l_{i,k}}}\bigg)^2
\end{align*}
we have the estimate:
\begin{equation}\label{e.areaest}
\snorm{A_{u,v}}^{p/2} \leq 2^{(p/2)-1} \bigg[ \bigg(\sum_{k=1}^{\infty} \sum_{i=1,2} \snorm{A_{l_{i,k},r_{i,k}}}\bigg)^{p/2} + \half \bigg(\sum_{k=1}^{\infty} \sum_{i=1,2} \snorm{X_{r_{i,k}}-X_{l_{i,k}}}\bigg)^{p} \bigg].
\end{equation}
Using H\"older's inequality, with $p>2$ and $\gamma>p-1$, we have
\begin{align}
\snorm{A_{u,v}}^{p/2} &\leq 2^{(p/2)-1} \bigg[ \bigg( \sum_{n=1}^{\infty} n^{-\gamma/((p/2)-1)} \bigg)^{(p/2)-1} \; \sum_{n=1}^{\infty} n^{\gamma} \bigg( \sum_{i=1,2} \snorm{A_{l_{i,k},r_{i,k}}} \bigg)^{p/2}\nonumber\\
&\qquad \qquad + \half \;\bigg( \sum_{n=1}^{\infty} n^{-\gamma/(p-1)} \bigg)^{p-1} \; \sum_{n=1}^{\infty} n^{\gamma} \bigg( \sum_{i=1,2} \snorm{X_{r_{i,k}}-X_{l_{i,k}}} \bigg)^{p}  \bigg]\nonumber\\
&\leq C_1(p,\gamma)\; \sum_{n=1}^{\infty} n^{\gamma} \sum_{i=1,2} \snorm{A_{l_{i,k},r_{i,k}}}^{p/2} \nonumber\\
&\qquad+ C_2(p,\gamma) \;  \sum_{n=1}^{\infty} n^{\gamma} \sum_{i=1,2} \snorm{X_{r_{i,k}}-X_{l_{i,k}}}^{p}\label{e.est}.
\end{align}

One can uniformly bound $\snorm{A_{u,v}}^{p/2}$ for any pair of times $u<v \in [0,T]$ by extending the estimate in \eqref{e.est} over all the dyadic intervals at each level $n$, that is,
\begin{align*}
 \snorm{A_{u,v}}^{p/2}&\leq  C_1(p,\gamma)\; \sum_{n=1}^{\infty} n^{\gamma} \sum_{i=1}^{2^n} \snorm{A_{l_{i,k},r_{i,k}}}^{p/2}\\
&\qquad\qquad + C_2(p,\gamma) \;  \sum_{n=1}^{\infty} n^{\gamma} \sum_{i=1}^{2^n} \snorm{X_{r_{i,k}}-X_{l_{i,k}}}^{p}.
\end{align*}
If the right hand side is finite a.s. then the area can be defined for any pair of times.

The $(p/2)$-variation of the L\'evy area can be estimated by the same bound.
\begin{align}\label{e.est2}
\sup_{\pi} \sum_{\pi} \snorm{A_{u,v}}^{p/2}&\leq  C_1(p,\gamma)\; \sum_{n=1}^{\infty} n^{\gamma} \sum_{i=1}^{2^n} \snorm{A_{l_{i,k},r_{i,k}}}^{p/2}\nonumber\\
&\qquad\qquad + C_2(p,\gamma) \;  \sum_{n=1}^{\infty} n^{\gamma} \sum_{i=1}^{2^n} \snorm{X_{r_{i,k}}-X_{l_{i,k}}}^{p}.
\end{align}

We use \eqref{e.2ndmom} to control the first sum
\begin{equation*}
\expn{\snorm{A_{s,t}}^{p/2}} \leq C\; (t-s)^{p/2}  \qquad\text{for}\; p\leq 4.
\end{equation*}
So we have
\begin{align*}
\expn{\sum_{n=1}^{\infty} n^{\gamma} \sum_{i=1}^{2^n} \snorm{A_{l_{i,k},r_{i,k}}}^{p/2}} &\leq C \sum_{n=1}^{\infty}n^{\gamma} \sum_{i=1}^{2^n} (2^{-n}T)^{p/2}\\
&= C \sum_{n=1}^{\infty}n^{\gamma} 2^{-n((p/2)-1)}\\
&< \infty \qquad \text{for}\; p>2.
\end{align*}
This implies that the first term in the right hand side of \eqref{e.est2} is a.s. finite.  Now we consider the second term of \eqref{e.est2}.

\begin{lem}\label{l.finsum}
\begin{equation*}
\sum_{n=1}^{\infty} n^{\gamma} \sum_{k=0}^{2^n -1} \snorm{X_{(k+1)2^{-n}T}-{X_{k2^{-n}T}}}^p\quad<\infty \quad\as
\end{equation*}
\end{lem}
Before proving the lemma we recall a result of Monroe \cite{mon2}.
\begin{defn}\label{df.min}
Let $B_t$ be a Brownian motion defined on a probability space $(\Omega, \sigf, \P)$.  A stopping time $T$ is said to be minimal if for any stopping time $S\leq T$, $B(T) \indist B(S)$ implies that a.s. $S=T$.
\end{defn}

\begin{thm}\cite[Theorem 11]{mon2}\label{thm.monroe11}
Let $(M_t)_{t\geq0}$ be a right continuous martingale.  Then there is a Brownian motion $(\Omega, \siggt, B_t)$ and a family $(T_t)$ of $\siggt$-stopping times such that the process $B_{T_t}$ has the same finite distributions as $M_t$.  The family $T_t$ is right continuous, increasing, and for each $t$, $T_t$ is minimal.  Moreover, if $M_t$ has stationary independent increments then so does $T_t$.
\end{thm}
\begin{rem}
It should be noted that the stopping times $T_t$ are not generally independent of $B_t$.  However, in the case of $\alpha$-stable processes $0< \alpha<2$ one can use subordination to gain independence of the stopping times \cite{boch}.
\end{rem}

\proofof{Lemma \ref{l.finsum}}
Let $(\tau_t)_{t\geq0}$ denote the collection of minimal stopping times for which
\begin{equation*}
X_t \indist B_{\tau_t}.
\end{equation*} 
The proof will be completed once it has been shown that 
\begin{equation}\label{e.lem1}
\sum_{n=1}^{\infty} n^{\gamma} \sum_{k=0}^{2^n -1} \snorm{B_{\tau((k+1)2^{-n}T)}-{B_{\tau(k2^{-n}T)}}}^p\quad<\infty \quad\as
\end{equation}
The following inequality holds because Brownian motion is $(1/p^{\prime})$-H\"older continuous a.s. for $p^{\prime}>2$:  
\begin{align}\label{e.lem2}
\snorm{B_{\tau((t_{k+1,n})}-{B_{\tau(t_{k,n})}}}^p &\leq C \,\snorm{\tau(t_{k+1,n}) - \tau(t_{k,n})}^{{{p}\over{p^{\prime}}}}\\
& \qquad \forall\,k=0,\ldots,2^n-1,\forall\,n\geq 1, \;\as\nonumber
\end{align}
where $t_{k,n} \dfn k2^{-n}T$ and $2< p^{\prime} <p$. 

\cite[Theorem 1]{mon} shows that the index of the process $\tau(s)$ is half that of the L\'evy process.  Therefore, with probability one, $\tau(s)$ has finite $(1+\delta)$-variation for all $\delta>0$.  

\begin{thm}\cite[Theorem 5]{mon2}\label{thm.inf}
If $\tau$ is a minimal stopping time and $\E (B_{\tau}) =0$, then $\E(\tau) = \E(B_{\tau}^2)$. 
\end{thm}
Consequently the process $(\tau_t)_{t \geq 0}$ can be controlled in the following way:
\begin{equation}\label{e.finmom}
\expn{\tau_t} = \expn{B_{\tau_t}^2} = \expn{X_t^2} = t\,\int_{\snorm{x}<1} \snorm{x}^2 \;\nu(dx)
\end{equation}
where $\nu$ is the L\'evy measure corresponding to the process $X_t$.  From \eqref{e.finmom} and Theorem \ref{thm.monroe11} we note that the process $\tau_t$ is a L\'evy process whose L\'evy measure, say $\mu$, satisfies the following: 
\begin{equation*}
\int_0^1 \,x \;\mu(dx) < \infty.
\end{equation*}
From this result we deduce that the process $\tau_t$ a.s. has bounded variation.  From \cite[Theorem 5]{shtat} we note that there is a positive constant $A$ such that 
\begin{equation*}
\prob{\tau_t \leq A\;t \;,\; \forall t\geq0} = 1.
\end{equation*}
From the above bound and using the fact that $\tau$ has stationary independent increments one can show  
\begin{align*}
\prob{\tau(t_{k+1,n}) - \tau(t_{k,n}) \leq A (t_{k+1,n} - t_{k,n}) = A 2^{-n}\: \big\vert\;\tau(t_{k,n}) } &=1, \\
\prob{\bigcap_{n\geq1} \bigcap_{k\geq0}^{2^n-1} \bigg(\snorm{\tau(t_{k+1,n}) - \tau(t_{k,n})} \leq A \,2^{-n} \bigg)} &= 1.
\end{align*}
Returning to \eqref{e.lem2} we see that
\begin{align*}
\snorm{B_{\tau((t_{k+1,n})}-{B_{\tau(t_{k,n})}}}^p &\leq C\,\snorm{\tau(t_{k+1,n}) - \tau(t_{k,n})}^{{{p}\over{p^{\prime}}}}\\
&\leq C\, A\, 2^{-n({{p}\over{p^{\prime}}})}\\
\intertext{which implies that}
\sum_{n=1}^{\infty} n^{\gamma} \sum_{k=0}^{2^n -1} \snorm{B_{\tau((k+1)2^{-n}T)}-{B_{\tau(k2^{-n}T)}}}^p &\leq C\, A\,\sum_{n=1}^{\infty} n^{\gamma}2^{-n({{p}\over{p^{\prime}}}-1)}< \infty
\end{align*} 
due to $p^{\prime}$ being chosen in the interval $(2,p)$.\endproof

This lemma concludes the proof that the bound in \eqref{e.est2} is finite, which shows that the area process a.s. has finite $(p/2)$-variation.\endproof

In this section we have proved that the area process exists and has finite $(p/2)$-variation when $(X_t)_{t\geq 0}$ has the form \eqref{e.form}.  To prove Theorems \ref{thm.areaconv}, \ref{thm.parea} we note that a general L\'evy process has the form
$$
X_t = at + B_t + L_t + \sum_{\substack{ 0\leq s< t \\ \snorm{\Delta X_s}\geq1}} \Delta X_s \qquad\as
$$ 
So, we need to add area corresponding to the drift vector and the jumps of size greater than one.  However, this part of the L\'evy process has bounded variation and is piecewise smooth so there is no problem defining its area.  Similarly, it has a.s. finite $(p/2)$-variation. 

\appendix
\section{Homeomorphic flows}\label{app.hom}

In this section we give a proof that the solutions, generated by \eqref{e.1} as the initial condition is varied, form a flow of homeomorphisms when the integrator is a continuous function.  The proof modifies the one given in \cite{tel5} for the existence and uniqueness of solution to \eqref{e.1}.  The main idea is that one uniformly bounds a sequence of iterated maps which have projections giving the convergence of the solutions with two different initial points and bounding the difference of the solutions.

First, we need some notation. 

\begin{defn}
Let $T^{(n)}(\R^d)$ denote the truncated tensor algebra of length $n$ over $\R^d$.  That is
\begin{equation*}
T^{(n)}(\R^d) \dfn \bigoplus_{i=0}^n (\R^d)^{\otimes i}
\end{equation*}
where $(\R^d)^{\otimes 0} = \R$ and $T^{(\infty)}(\R^d)$ denotes the tensor algebra over $\R^d$.

Let $\Delta = [0,T]\times [0,T]$.  A map $X:\Delta \rightarrow T^{(n)}(\R^d)$ will be called a multiplicative functional of size $n$ if for all times $s<t<u$ in $[0,T]$ the following relation holds in  $T^{(n)}(\R^d)\,:$
\begin{equation*}
X_{st}\otimes X_{tu} = X_{su}
\end{equation*}
and $X^{(0)}_{st} \equiv 1$.

A map $X : \Delta \rightarrow T^{(n)}(\R^d)$ is called a classical multiplicative functional if $t\rightarrow X_t \dfn X^{(1)}_{0t}$ is continuous and piecewise smooth and 
\begin{equation}\label{e.class}
X^{(i)}_{st} = \iint_{s<u_1<\ldots<u_i<t} \;dX_{u_1}\ldots dX_{u_i}
\end{equation}
where the right hand side is a Lebesgue-Stieltjes integral.  We denote the set of all classical multiplicative functionals in $T^{(n)}(\R^d)$ by $S^{(n)}(\R^d)$. 
\end{defn}

\begin{defn}\label{d.control}
We call a continuous function $\omega:\Delta \rightarrow \R^+$ a control function if it is super-additive and regular, that is,
\begin{alignat*}{3}
\omega\,(s,t) + \omega\,(t,u) &\leq \omega\,(s,u)& \qquad\forall\; s<t<u\in [0,T]\\
\omega\,(s,s) &= 0 &\forall \;s\in [0,T]\,.
\end{alignat*}
\end{defn}
\begin{eg}
Let $X$ be a path of strong finite $p$-variation.  Then we can define the following control function:
\begin{equation}
\ost \dfn \variationt{X}{p}{s}{t}^p.
\end{equation}
\end{eg}
\begin{defn}\label{d.enhp}
A functional $X = (1,X^{(1)},\ldots,X^{(n)})$ defined on $T^{(n)}(\R^d)$ where $n=[p]$ is said to have finite $p$-variation if there is a control function $\omega$ such that 
\begin{equation}\label{e.varcon}
\snorm{X_{st}^{(i)}} \leq {{\ost^{i/p}}\over{\beta (i/p)!}} \qquad \forall\; (s,t)\in \Delta,\; i=1,\ldots,n
\end{equation}
for some sufficiently large $\beta$ and  $x! \dfn \Gamma(x+1)$.  
\end{defn}
\begin{thm}\cite[Theorem 2.2.1]{tel5}\label{thm.enhance}
Let $X^{(n)}$ be a multiplicative functional of degree $n$ which has finite $p$-variation, with $n\dfn [p]$ ($[p]$ denotes the integer part of $p$).  Then for $m>n$ there is a unique multiplicative extension $X^{(m)}$ in $T^{(m)}(\R^d)$ which has finite $p$-variation.
\end{thm}

\begin{rem}
The above theorem shows that once a sufficient number of low order integrals associated to a path $X_t$ have been defined, then the remaining iterated integrals of $X_t$ are defined.  
\end{rem}

\begin{defn}
We call a multiplicative functional $ X:\Delta\rightarrow T^{(n)}(\R^d)$ geometric if there is a control function $\omega$ such that for any positive $\epsilon$ there exists a classical multiplicative functional $ Y(\epsilon)$ which approximates $X$ in the following way:
\begin{equation*}
\sllnorm{\bigg( X_{st}- Y_{st}(\epsilon)\bigg)^{(i)}} \leq \epsilon \,\ost^{i/p} \qquad i=1,\ldots,n=[p]\,.
\end{equation*}
We denote the class of geometric multiplicative functionals with finite $p$\linebreak[4]-variation by $\Omega G (\R^d)^p$.
\end{defn}
\begin{eg}
Let $W_t$ be an $\R^d$-valued Brownian motion.  Then the following functional $ W$ defined on $T^{(2)}(\R^d)$ belongs to $\Omega G (\R^d)^p$ for any $p>2$.
\begin{equation}\label{e.bmgmf}
 W_{st} \dfn \bigg(1,W_t - W_s, \iint_{s<u_1<u_2<t} \circ dW_{u_1}\, \circ dW_{u_2} \bigg)
\end{equation}
where $\circ dW_u$ denotes the Stratonovich integral.  It should be noted that if one replaced the Stratonovich differential in \eqref{e.bmgmf} by the It\^o differential then one would not get an element of $\Omega G (\R^d)^p$.  This is due to the quadratic variation term which occurs in the symmetric part of the area process $$ W^{(2)}_{st}=\iint_{s<u_1<u_2<t}  dW_{u_1}\,  dW_{u_2}.$$  
\end{eg}

It was shown in \cite{sip} that one had sufficient control of the above functional to generate path-wise solutions to SDEs driven by a Brownian motion.  This control was derived from a moment condition in the same spirit as Kolmogorov's criterion for H\"older continuous paths.  The moment condition was verified for the above area by the use of known stochastic integral results, though one could also derive it from a construction depending on the linearly interpolated Brownian motion.  

There are two stages to defining the integral against a geometric multiplicative functional.  The first gives a functional which is almost multiplicative (see \cite{tel5} for definition).  The second associates, uniquely, a multiplicative functional to the almost multiplicative functional. 

\begin{thm}\cite{tel5}
There is a unique geometric multiplicative functional $Y$ which we call the integral of the 1-form $\theta$ against the geometric multiplicative functional $ X$.  We denote this by
\begin{equation*}
Y_{st} \dfn \int_s^t \theta (X_u) \;\delta  X.
\end{equation*}
\end{thm}

\begin{cor}\label{cor.pint}
One has the following control on the $p$-variation of $Y$:
\begin{equation}
\sllnorm{\bigg(\int_s^t \theta (X_u) \;\delta  X\bigg)^{(i)}} \leq \bigg(C \;\ost\bigg)^{i/p}/\bigg(\beta \bigg(i/p\bigg)!\bigg)\qquad i=1,\ldots,[p]
\end{equation}
where $C$ depends on $p, \norm{f}{\lip(\gamma)}, \gamma, \lambda, \beta,L$ and $[p]$.
\end{cor}

The estimate is derived from estimating both the almost multiplicative functional and the difference of it from the integral. 

We now state two lemmas which help prove that the solutions of \eqref{e.1} are homeomorphic flows when the initial condition is varied. 

\begin{lem}\label{lem.gron}
Let $X$ be in $\geo{\R^d}$ controlled by a regular $\omega_0$.  Let $f:\R^n\rightarrow \hom(\R^d,\R^n)$ be a $\lip(\gamma)$ map for some $\gamma>p$.  Let $ Y^{(i)}_{st},\, i=1,2$ denote the element in $\geo{\R^n}$ which solves the rough integral equation
$$
Y^{(i)}_{st} = \int_s^t f(Y^{(i)})\,\delta X
$$
with initial condition $Y^{(i)}_{0}=a_i,\, i=1,2$.
Let $ W_{st}$ be the multiplicative functional which records the difference in the multiplicative functionals $ Y^{(1)}_{st}$ and $ Y^{(2)}_{st}$. Then
\begin{equation}\label{e.origbd}
\slnorm{ W^{(i)}_{st}} \leq   \theta^i  {{\ost^{(i/p)}}\over{\beta(i/p)!}}\qquad \forall\; i \geq 1,
\end{equation}
where $\theta = \snorm{a_1 -a_2}$,  $\omega \dfn C\,\omega_0$, the constant $C$ depends on $p,$ $\norm{f}{\lip(\gamma)},$ $\beta,$ $\gamma$.  The bound holds for all times $s\leq t$ on the interval $J\dfn \{u\,:\,\omega\,(0,u) \leq 1\}$.
\end{lem}

\begin{lem}\label{lem.gron2}
With the assumptions of Lemma \ref{lem.gron} one can estimate the difference of the increments of $Y^{(1)}_{st}$ and $Y^{(2)}_{st}$ for any pair of times $0\leq s<t$ which satisfy $\ost \leq 1$ as follows:
\begin{align*}
\snorm{Y^{(1)}_{st}-Y^{(2)}_{st}} &\leq \theta\,\exp\bigg[ {{1}\over{\beta(1/p)!}}\,\big( \omega\,(0,s) + \omega\,(0,s)^{(1/p)}\big)\bigg]\;{{\ost^{(1/p)}}\over{\beta(1/p)!}}
\end{align*}
In particular for any $t>0$ one has:
\begin{equation}\label{e.gro}
\snorm{Y^{(1)}_t - Y^{(2)}_t} \leq \snorm{a_1 -a_2}\; C(t).
\end{equation} 
\end{lem}

Now we can prove that the solutions form a flow of homeomorphisms as the initial condition is varied.

\proofof{Theorem \ref{thm.hom}}
The continuity of solutions follows from Lemma \ref{lem.gron2}.  It remains to show that the inverse map exists and is continuous.  This can be checked by repeating all the previous arguments using the reversed path $(X_{t-s})_{0\leq s\leq t}$ as the integrator.\endproof

The induction part of the proof of Lemma \ref{lem.gron} will require the following lemma about rescaling:

\begin{lem}\label{lem.scal}\cite{tel5}
Let $ X$ be a multiplicative functional in $T^{([p])}(\R^d)$ which is of finite $p$-variation controlled by $\omega$.  Let $( X, Y)$ be an extension of $ X$ to $T^{([p])}(\R^d\oplus\R^n)$ of finite $p$-variation controlled by $K \omega$.  Then $( X, \phi  Y)$ is controlled by
\begin{equation*}
\max \bigg\{ 1, \phi^{kp/i} K : 1 \leq k \leq i \leq [p] \bigg\} \omega
\end{equation*}
where $\phi \in \R$.  In particular, if $\phi \leq K^{-[p]/p} \leq 1$ then $( X, \phi  Y)$ is controlled by $\omega$.
\end{lem}

\proofof{Lemma \ref{lem.gron}}
We set up an iteration scheme of multiplicative functionals which we will bound uniformly, by induction.  A projection of the sequence proves that a Picard iteration scheme converges to the solutions of \eqref{e.1} starting from $a_1$ and $a_2$.  Another projection shows that the difference of these solutions is bounded. 

Let $\epsilon>0$ and $\eta >1$.  Let $V^{(1)}_{st}$ be the geometric multiplicative functional given by
\begin{align*}
V^{(1)}_{st}&\dfn(Z_{st}^{(1)(1)},\,Y_{st}^{(1)(1)},\,Y_{st}^{(1)(0)},\,Z_{st}^{(2)(1)},\,Y_{st}^{(2)(1)},\,Y_{st}^{(2)(0)},\,W_{st}^{(1)},\,\epsilon^{-1}X_{st})\\
&=\bigg(\int_s^t\,f(a_1)\;\delta X -a_1,\,\int_s^t\,f(a_1)\;\delta X,\,a_1,\,\int_s^t\,f(a_2)\;\delta X -a_2,\,\\
&\qquad\int_s^t\,f(a_2)\;\delta X,\,a_2,\, \int_s^t\,f(a_1)-f(a_2)\;\delta X,\,\epsilon^{-1}X_{st}\bigg).
\end{align*}
The iteration step is a two stage process.  Given $V^{(m)}$ we set 
$$
\tilde V^{(m+1)} = \int k^{m}_{\theta} (V^{(m)})\;\delta V^{(m)}
$$
where $k^m_{\theta}$ is the $1$-form on $((\R^n)^{\oplus7}\oplus \R^d)$ given by
\begin{align*}
&k^m_{\theta}(a_1,\ldots,a_8)\,(dA_1,\ldots,dA_8)\\
&\qquad = \bigg( a_1\,g(a_2,a_3)\,dA_8,\, dA_3 + \eta^{-m} dA_1,\, dA_2,\,a_4\,g(a_5,a_6)\,dA_8,\\
&\qquad\qquad\qquad \, dA_6 + \eta^{-m} dA_4, \, dA_5, \, \theta^{-1}\,g(a_2,a_4)\,dA_8,\, dA_8 \bigg).
\end{align*}
$g(x,y)$ is the 1-form appearing in \cite[Lemma 3.2]{tel1} which satisfies the following relation with respect to $f$:
\begin{equation*}\label{e.g}
f^i(x)-f^i(y) = \sum_j (x-y)^j\;g^{ij}(x,y).
\end{equation*}
$\tilde V^{(m+1)}$ is well defined because $g$ and $k^m_{\theta}$ are both $\lip(\gamma)$ for some $\gamma>p-1$. 

We define $V^{(m+1)}$ to be the geometric multiplicative functional obtained by rescaling the first and fourth components of $\tilde V^{(m+1)}$ by $\epsilon\eta$ and the seventh component by $\epsilon$.

The uniform bound on the iterates $(V^{(m)})_{m\geq1}$ will be obtained by induction.  $X$ is controlled by a regular $\omega_0$ so there exists a constant $C$ such that $V^{(1)}$ is controlled by $\omega \dfn C\,\omega_0$.  Suppose that $V^{(k)} (k\leq m)$ are controlled by $\omega$.  From (Corollary \ref{cor.pint}) there is a constant $C_1$ such that $\tilde V^{(m+1)}$ is controlled by $C_1 \omega$.  If we choose $\epsilon>0,\,\eta>1$ such that $\epsilon \leq C_1^{{{-[p]}\over{p}}}$ and $\epsilon \eta \leq C_1^{{{-[p]}\over{p}}}$, then Lemma \ref{lem.scal} implies that $V^{(m+1)}$ is controlled by $\omega$, completing the induction step.

The uniform control on the iterates $ V^{(m)}$ ensures the convergence of \linebreak $\{ Y^{(i)(m)}\}_{m\geq 1}$ to the solutions of 
\begin{align*}
dY^{(i)}_t &=f(Y^{(i)}_t)\;dX_t \qquad Y^{(i)}_0=a_i,\;i=1,2.
\end{align*}

Through the definition of $\{{}^{\theta} W^{(m)}\}_{m\geq 1}$, the sequence at the level of the paths will converge to the scaled difference of the two solutions $\theta^{-1} (Y^{(2)}-Y^{(1)})$.  For $s,t$ in $J$ one has 
\begin{equation*}
\snorm{{}^{\theta} W_{st}^{(i)}} \leq {{\ost^{i/p}}\over{\beta (i/p)!}}         \qquad i=1,\ldots,[p],
\end{equation*}
which implies that
  \begin{equation*}
\snorm{ W_{st}^{(i)}} \leq \theta^i {{\ost^{i/p}}\over{\beta (i/p)!}}         \qquad i=1,\ldots,[p].\eendproof
\end{equation*}

\proofof{Lemma \ref{lem.gron2}}
We define the following set of times:
\begin{align}\label{e.times}
t_0 \dfn 0 &\text{ and } t_j \dfn \inf\{ u > t_{j-1} \,:\, \omega\,(t_{j-1},u) = 1\} \qquad \forall\; j\in\{1,\ldots,n(s)\}\\
\text{where} \;&n(s) \dfn \max\{j\,:\, t_j \leq s\}\;\text{and}\; t_{n(s)+1}=s\nonumber\,.
\end{align}

We solve the differential equation starting from $s$ and use \eqref{e.origbd} to show that 
\begin{equation*}
\snorm{ W^k_{st}} \leq K(s)^k\,{{\ost^{(k/p)}}\over{\beta(k/p)!}}
\end{equation*}
where $K(s)$ is an upper bound on the supremum over all the possible differences of the paths $\snorm{ Y^{(1)}_s -  Y^{(2)}_s}\,,$ at time $s$.   The bound $K(s)$ is derived recursively by considering the analogous upper bound for the difference of the solutions to the differential equation over the time interval $[t_{i-1},t_i]$ given below:
\begin{align}
\snorm{Y^{(1)}_{t_{i}}-Y^{(2)}_{t_{i}}} &\leq \snorm{Y^{(1)}_{t_{i-1}}-Y^{(2)}_{t_{i-1}}} + \snorm{W_{t_{i-1}\,t_i}}\nonumber\\
&\leq \snorm{Y^{(1)}_{t_{i-1}}-Y^{(2)}_{t_{i-1}}}\bigg( 1 + {{\omega\,(t_{i-1},t_i)^{(1/p)}}\over{\beta(1/p)!}} \bigg)\nonumber\\
\intertext{which implies that}
K(t_j) &\leq K(t_{j-1})\;\bigg\{ 1+ {{\omega\,(t_{j-1},t_j)^{(1/p)}}\over{\beta(1/p)!}}\bigg\}\qquad j=1,\ldots,n(s)+1.\nonumber
\end{align}
Therefore
\begin{align*}
\snorm{ W^k_{st}} &\leq K(t_0)^k\,\prod_{j=1}^{n(s)+1} \, \bigg\{ 1+ {{\omega\,(t_{j-1},t_j)^{(1/p)}}\over{\beta(1/p)!}}\bigg\}^k\;{{\ost^{(k/p)}}\over{\beta(k/p)!}}\\
&\leq \theta^k\,\exp\bigg[ k\,\big(\sum_{j=1}^{n(s)} {{\omega\,(t_{j-1},t_j)^{(1/p)}}\over{\beta(1/p)!}} + {{\omega\,(t_{n(s)},s)^{(1/p)}}\over{\beta(1/p)!}}\big)\bigg]\;{{\ost^{(k/p)}}\over{\beta(k/p)!}},
\end{align*}
noting that $\omega(t_{j-1},t_j) = 1$ and using the sub-additivity of $\omega$ we obtain
\begin{align*}
&\leq \theta^k\,\exp\bigg[ {{k}\over{\beta(1/p)!}}\,\big( \omega\,(0,s) + \omega\,(0,s)^{(1/p)}\big)\bigg]\;{{\ost^{(k/p)}}\over{\beta(k/p)!}}.
\end{align*}
By considering the above bound at the level of the paths $(k=1)$ and repeatedly using the triangle inequality one deduces \eqref{e.gro}.\endproof

\bibliographystyle{plain}
\bibliography{journal,institut,publish,ref}
\address
\end{document}